\documentclass[11pt]{amsart}
\usepackage{geometry}                % See geometry.pdf to learn the layout options. There are lots.
\usepackage{graphicx}
\usepackage{url}
\usepackage{amssymb}
\usepackage{amsrefs}
\usepackage{epstopdf}
\usepackage{tikz}
\usetikzlibrary{shapes,patterns,calc,snakes}

\newcommand\mbb{\mathbb}

\newcommand\C{\mbb{C}}

\newcommand\N{\mbb{N}}

\newcommand\R{\mbb{R}}

\theoremstyle{plain}
\newtheorem{Thm}[equation]{Theorem}

\newtheorem*{Thm*}{Theorem}
\newtheorem*{Prop*}{Proposition}
\newtheorem*{Cor*}{Corollary}
\newtheorem*{Lemma*}{Lemma}
\newtheorem*{Sublemma*}{Sublemma}
\newtheorem*{Conjecture*}{Conjecture}

\theoremstyle{definition}

\newtheorem*{Construction*}{Construction}
\newtheorem*{Def*}{Definition}
\newtheorem*{Defs*}{Definitions}
\newtheorem*{Example*}{Example}
\newtheorem*{Examples*}{Examples}
\newtheorem*{LemmaDef*}{Lemma and Definition}
\newtheorem*{Notation*}{Notation}
\newtheorem*{Problem*}{Problem}
\newtheorem*{Question*}{Question}
\newtheorem*{Remark*}{Remark}
\newtheorem*{Remarks*}{Remarks}
\newtheorem*{Warning*}{Warning}

\title{Free Semialgebraic Geometry}
\author{Tim Netzer}
\begin{document}

\begin{abstract}
This is a survey article on the currently very active research area of free (=non-commutative) real algebra and geometry. We first review some of the important results from the commutative theory, and then explain similarities and differences as well as some important results in the free setup.
\end{abstract}

\maketitle
\section{Classical Semialgebraic Geometry}
In this section we briefly review some concepts and results from commutative real algebra and geometry. For details and proofs see for example \cite{pd, m, bcr, netznote, habil}.

Important objects of study in classical (=commutative) real algebra and geometry are semialgebraic sets. A {\it basic closed semialgebraic set} is of the form $$W(p_1,\ldots, p_r):=\left\{ a\in\R^d\mid p_1(a)\geq 0, \ldots, p_r(a)\geq 0\right\}$$ where $p_1,\ldots, p_r\in\R[x_1,\ldots, x_d]$ are polynomials.  A general {\it semialgebraic set} is a (finite) Boolean combination of basic closed semialgebraic sets. 
\begin{figure}[h!]
\begin{tikzpicture}
    \begin{scope}
    \node {\includegraphics[width=1 in]{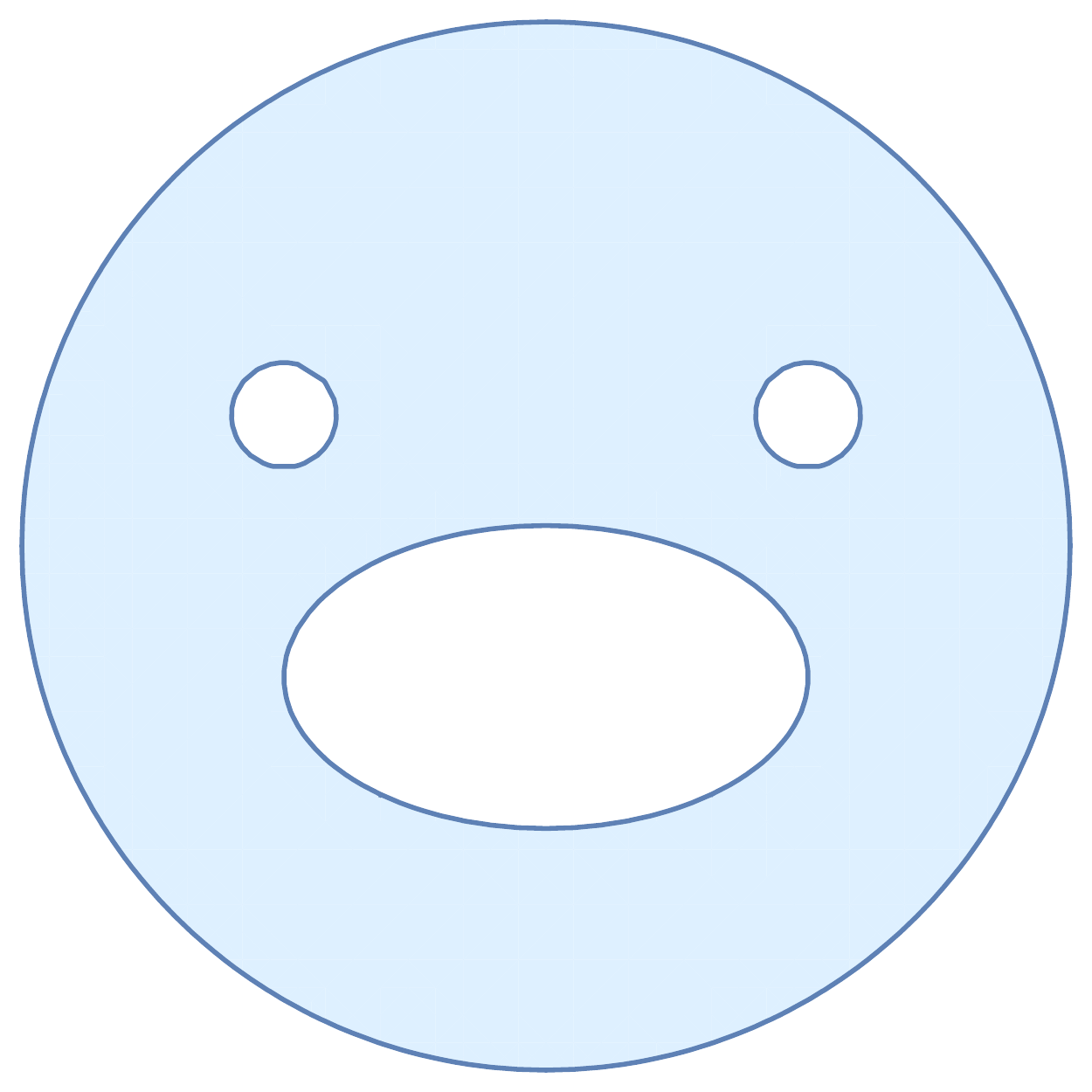}};
    \end{scope}
    \begin{scope}[xshift=5cm]
    \node {\includegraphics[width=1.3 in]{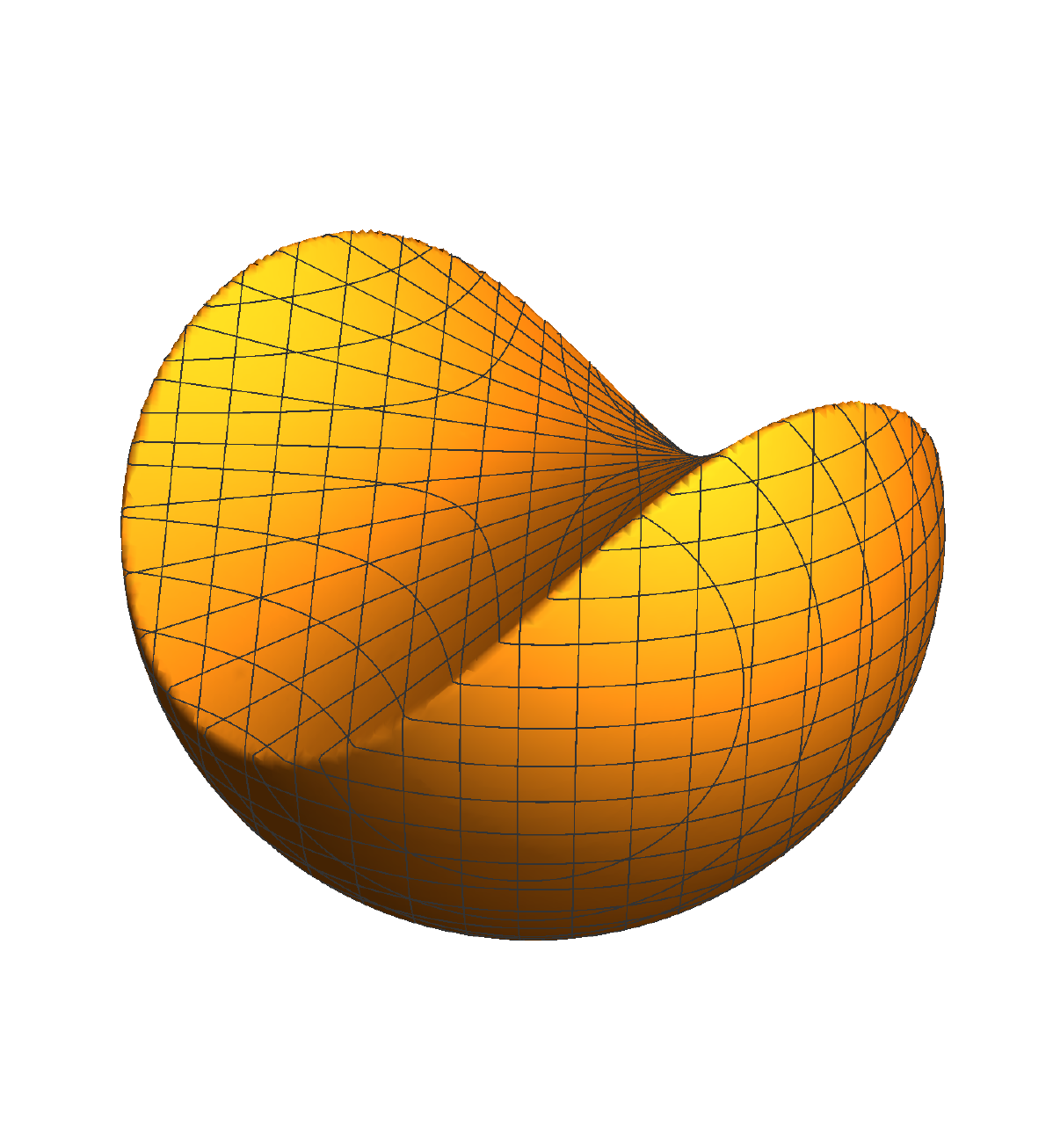}};
    \end{scope}
\end{tikzpicture}
\caption{Two (basic closed) semialgebraic sets}
\end{figure}
An important and foundational result in real algebraic geometry is the following:
\begin{Thm}[Projection Theorem] Projections of semialgebraic sets are again semialgebraic.
\end{Thm}
The Projection Theorem is not easy to proof. It also has some strong implications for logic and model theory of real closed fields. Since projections correspond to existential quantifiers, it leads to quantifier elimination in the theory of real closed fields, which lies at the core of almost any Positivstellensatz in real algebra.  It also proves decidability of the theory of real closed fields. Although not easy, the Projection Theorem admits constructive proofs. In practice, finding a semialgebraic description of a projection can be very challenging, however.  

In classical algebraic geometry, affine varieties are classified via polynomial functions that vanish on them. These functions are described algebraically by Hilbert's Nullstellensatz.  For semialgebraic sets, one considers  {\it nonnegative} polynomials, and  {\it Positivstellens\"atze}  provide algebraic characterization of such polynomials. The role of ideals is  taken by {\it preorderings}. The preordering generated by  polynomials $p_1,\ldots, p_r\in\R[x_1,\ldots, x_d]$ arises from the $p_i$ and sums of squares of polynomials, by addition and multiplication:

$${\mathcal P}(p_1,\ldots, p_r):=\left\{ \sum_{e\in\{0,1\}^r} \sigma_e p_1^{e_1}\cdots p_r^{e_r} \mid \sigma_e\ \mbox{ sums of squares of polynomials}\right\}.$$ Note that polynomials from the preordering are obviously nonnegative on $W(p_1,\ldots, p_r)$.

\begin{Thm}[Nichtnegativstellensatz]
For $p,p_1,\ldots, p_r\in\R[x_1,\ldots, x_d]$, the following are equivalent:
\begin{itemize}
\item[($i$)] $p\geq 0$ on $W(p_1,\ldots, p_r)$.
\item[($ii$)] $t_1p=p^{2n}+t_2$ for some $t_1,t_2\in {\mathcal P}(p_1,\ldots, p_r)$, $n\in\N$.
\end{itemize}
\end{Thm} 

The direction ($ii$)$\Rightarrow$($i$) is straightforward to see, so ($ii$) is an algebraic certificate for nonnegativity of $p$. The factor $t_1$ in ($ii$) is often called a {\it denominator}. Over the field $\mathbb R(x_1,\ldots, x_d)$ it can be brought to the other side, providing a representation with rational functions.  The case $r=1$ and  $p_1=1$ is precisely Hilbert's 17th Problem: every globally nonnegative polynomial is a sum of squares of rational functions.
One can get rid of the denominator only under additional assumptions. The first and most important such Positivstellensatz without denominators is the following, where the conditions of boundedness and strict positivity is necessary for the theorem to hold:
\begin{Thm}[Schm\"udgen's Positivstellensatz]
Let $p_1,\ldots, p_r\in\R[x_1,\ldots, x_d]$ and assume $W(p_1,\ldots,p_r)$ is bounded. Then for any $p\in \R[x_1,\ldots, x_d]$ $$p>0 \mbox{ on } W(p_1,\ldots, p_r) \quad\Rightarrow\quad p\in  {\mathcal P}(p_1,\ldots, p_r).$$
\end{Thm}
After this very brief introduction, let us now pass to the non-commutative setup.

\section{Free Real Algebra and Geometry}

Semialgebraic sets are defined by polynomial inequalities. So before we can talk about non-commutative semialgebraic sets, we introduce non-commutative polynomials. In the non-commutative setup we will always use an involution (in fact the involution is also there in the classical case, however invisible since it is just the identity on real polynomials). In presence of an involution we can use complex numbers as our ground field and take Hermitian elements as the "real"\ objects. Using complex numbers is often more convenient and allows for cleaner proofs.

So let $\C\langle z_1,\ldots, z_d\rangle$ denote the algebra of polynomials in the {\it non-commuting variables} $z_1,\ldots, z_d$. Their elements are $\C$-linear combinations of {\it words} in the variables. Since the variables do not commute, words like $z_1z_2$ and $z_2z_1$ are different, where they would coincide in the commutative case. We use the involution $*$ on $\C\langle z_1,\ldots, z_d\rangle$ that fixes the variables, i.e.\ $z_i^*=z_i$ holds for all $i$, but reverses the order in each word and acts as complex conjugation on the coefficients. For example, we have $$\left(z_1^3-z_1z_2+i\right)^*=z_1^3-z_2z_1-i.$$
Let $$\C\langle z_1,\ldots, z_d\rangle_h:=\left\{ p\in \C\langle z_1,\ldots, z_d\rangle\mid p^*=p\right\}$$ be the set of {\it Hermitian elements}. They form a real vectorspace, but not an algebra (in case $d\geq 2)$. Note that Hermitian elements do not necessarily have real coefficients, and polynomials with real coefficients are not necessarily Hermitian.

Into a non-commutative polynomial $p\in \C\langle z_1,\ldots, z_d\rangle$ we can plug in a $d$-tuple of elements from any complex algebra, and obtain an element from this algebra as the result. We will restrict ourselves to matrix algebras here, i.e. we take $(A_1,\ldots, A_d)\in {\rm Mat}_s(\C)^d$ for some $s\geq 1$ and obtain $$p(A_1,\ldots, A_d)\in {\rm Mat}_s(\C).$$ The need for an involution and Hermitian elements becomes clear when trying to  capture real phenomena. Indeed if $A_1,\ldots,A_d\in {\rm Her}_s(\C)$ are Hermitian matrices and $p\in \C\langle z_1,\ldots, z_d\rangle_h$ is Hermitian as well, then so is the result: $$p(A_1,\ldots, A_d)\in{\rm Her}_s(\C).$$ A Hermitian matrix  is {\it positive semidefinite} if all of its Eigenvalues are nonnegative; we denote this by $\geqslant 0$. This is the right notion of positivity in our setup, so if $$p(A_1,\ldots, A_d)\geqslant 0$$ we say that $p$ is nonnegative at the (non-commutative) point $(A_1,\ldots, A_d)\in{\rm Her}_s(\C)^d$.

 It is obvious that every matrix  $A\in{\rm Her}_s(\C)$ that can be written as an {\it Hermitian square} $$A=B^*B$$ for some $B\in {\rm Mat}_s(\C)$ is positive semidefinite; in fact every positive semidefinite matrix is of that form. So if $$p=\sum_{i=1}^nq_i^*q_i$$ for certain $q_1,\ldots, q_n\in  \C\langle z_1,\ldots, z_d\rangle,$ we obtain $$p(A_1,\ldots, A_d)=\sum_i q_i(A_1,\ldots, A_d)^*q_i(A_1,\ldots, A_d)\geqslant 0$$ for any  $(A_1,\ldots, A_d)\in{\rm Her}_s(\C)^d$. So the set of {\it sums of Hermitian squares} $$\Sigma\C\langle z_1,\ldots, z_d\rangle^2:=\left\{ \sum_{i=1}^n q_i^*q_i\mid n\in\N, q_i\in \C\langle z_1,\ldots, z_d\rangle\right\} \subseteq \C\langle z_1,\ldots, z_d\rangle_h$$ only contains polynomials that are positive semidefinite on each Hermitian matrix tuple. The first surprising result, a global Positivstellensatz and a non-commutative analogue of Hilbert's 17th Problem, is due to Helton:

\begin{Thm}[\cite{he0}]\label{hth} Let $p\in\C\langle z_1,\ldots, z_d\rangle_h$ and assume $$p(A_1,\ldots, A_d)\geqslant 0$$ for all $(A_1,\ldots, A_d)\in{\rm Her}_s(\C)^d$ and all $s\geq 1$. Then $$p\in\Sigma\C\langle z_1,\ldots, z_d\rangle^2.$$
\end{Thm}

In contrast to the commutative result, no denominator is needed in Helton's theorem. However, the natural notion of positivity is  much stronger here than in Hilbert's 17th Problem, where positivity is only assumed on matrices of size $1$, instead of matrices of all sizes. Note however that also in Helton's result one can bound the matrix size,  depending only on $d$ and the degree of $p$.

Let us now define free basic closed semialgebraic sets. In analogy to the above described commutative setup, we define for $p_1,\ldots, p_r\in\C\langle z_1,\ldots, z_d\rangle_h$ and $s\geq 1$ $$W_s(p_1,\ldots, p_r):=\left\{ (A_1,\ldots, A_d)\in{\rm Her}_s(\C)^d\mid p_i(A_1,\ldots, A_d)\geqslant 0, i=1,\ldots, r\right\}.$$ A guiding principle in non-commutative geometry is to not consider matrices of one size alone, but all sizes at once. We thus define the {\it free basic closed semialgebraic set} defined by $p_1,\ldots, p_r$ as $${\rm F}W(p_1,\ldots, p_r):=\left(W_s(p_1,\ldots, p_r)\right)_{s=1}^\infty.$$
There is no known Positivstellensatz for positivity on general free basic closed semialgebraic sets, however certain results in  special cases. One of them deals with the {\it matrix cube}, see for example \cite{aleknetzthom}:

\begin{Thm}\label{box}
Assume  $p\in\C\langle z_1,\ldots, z_d\rangle_h$ is nonnegative on $${\rm F}W\left(1-z_1^2,\ldots, 1-z_d^2\right).$$ Then there exists a representation $$p=\sum_{i} q_i^*q_i + \sum_{i,j}q_{ij}^*\left(1-z_j^2\right)q_{ij}$$ for certain $q_i,q_{ij}\in\C\langle z_1,\ldots, z_d\rangle_h.$
\end{Thm}
Another such Positivstellensatz is explained in the next section, and we also refer to \cite{aleknetzthom} for more examples and unified proofs.

A  notion of free semialgebraic sets beyond free basic closed semialgebraic sets has not been established in the literature so far. Boolean combination of basic closed sets will surely have to be allowed, but maybe that is not yet the best possible notion. This becomes clear when trying to prove a free projection theorem. For any $s\geq 1$ we apply the projection map \begin{align*}\pi_s\colon {\rm Her}_s(\C)^d&\to{\rm Her}_s(\C)^{d-1} \\ (A_1,\ldots, A_d)&\mapsto (A_1,\ldots, A_{d-1})\end{align*} to $W_s(p_1,\ldots, p_r)$ and altogether obtain $$\pi({\rm F}W(p_1,\ldots, p_r)):=\left(\pi_s(W_s(p_1,\ldots, p_r))\right)_{s=1}^\infty.$$How does such a projected  set look like,  is it a Boolean combination of free basic closed sets? The answer to this question is no,  and the whole topic does not look too encouraging. For example (see \cite{dreschnetzthom}), using free basic closed semialgebraic sets, Boolean combinations and projections, one can construct the set $$\left(S_s\right)_{s=1}^\infty, \qquad S_s=\left\{\begin{array}{cc} {\rm Her}_s(\C) & s \mbox{ prime} \\ \emptyset & \mbox{ else.} \end{array} \right.$$ This set cannot be defined without projections, even if the language is enlarged by using trace, determinant and many other functions. Even more discouraging  is the following result from \cite{dreschnetzthom}:

\begin{Thm}
It is undecidable whether a set constructed from free basic closed semialgebraic sets, Boolean combinations and projections is empty (at each level).
\end{Thm}

A very recent positive result is \cite{kltr}. Without going too much into the details, it states that quantifiers in non-commutative formulas can be eliminated, as long as the formula is evaluated at matrix tuples of {\it fixed size}. This is not a trivial result, since the variables in such formulas do not refer to the single matrix entries (where the result would follow from classical (commutative) quantifier elimination), but to matrices as a whole. However, the formula without quantifiers will depend on the matrix size. So the result does  not imply a general (size independent) projection theorem.

A much more fruitful concept is {\it free convexity}, as we now explain in our last section.

\section{Free Convexity}
 
 In this section we define the notion of a non-commutative convex set. Since definitions and results become cleaner for convex cones instead of convex sets, we restrict ourselves to cones here. As above we consider free sets $$S=\left(S_s\right)_{s=1}^\infty, \quad S_s\subseteq {\rm Her}_s(\C)^d \ \mbox{ for all } s\geq 1.$$ Matrix convexity of $S$ is defined via two properties. A very reasonable assumption, even fulfilled for all free basic closed semialgebraic sets, is {\it closedness under direct sums}. For $\underline A=(A_1,\ldots, A_d)\in{\rm Her}_s(\C)^d, \underline B=(B_1,\ldots, B_d)\in{\rm Her}_t(\C)^d$ we define $$\underline A\oplus\underline B:=\left(\left(\begin{array}{c|c}A_1 & 0 \\\hline0 & B_1\end{array}\right), \ldots,\left(\begin{array}{c|c}A_d & 0 \\\hline0 & B_d\end{array}\right)\right) \in {\rm Her}_{s+t}(\C)^d.$$ $S$ is closed under direct sums if \begin{equation}\tag{C1}\underline A\in S_s, \underline B\in S_t\ \Rightarrow\ \underline A\oplus \underline B \in S_{s+t}.\end{equation} The second condition resembles scaling with positive reals, but even connects the different levels of $S$. For $V\in {\rm Mat}_{s,t}(\C)$ and $\underline A\in{\rm Her}_s(\C)^d$ we define $$V^*\underline AV:= \left(V^*A_1V,\ldots, V^*A_dV \right)\in {\rm Her}_t(\C)^d.$$ The second condition then reads \begin{equation}\tag{C2}\underline A\in S_s, V\in {\rm Mat}_{s,t}(\C)\ \Rightarrow\ V^*\underline A V\in S_t.\end{equation} If $S$ fulfills (C1) and (C2), it is called a {\it matrix convex cone}. It is easily checked that each $S_s$ is a classical convex cone in the real vector space ${\rm Her}_s(\C)^d$ in this case. However, matrix convexity is a stronger condition in general, connecting the different levels of $S$ via (C2). Also note that a matrix convex cone is almost the same as an  {\it abstract operator system} \cite{fnt, pau}, which only requires all $S_s$ to be closed and salient with nonempty interior, additionally.

The most basic examples of matrix convex cones are {\it free spectrahedral cones} (or  operator systems with a finite-dimensional realization, equivalently). For $M_1,\ldots, M_d\in {\rm Her}_r(\C)$ define $$S_s(M_1,\ldots, M_d):=\left\{ (A_1,\ldots, A_d)\in{\rm Her}_s(\C)^d\mid M_1\otimes A_1+\cdots +M_d\otimes A_d\geqslant 0\right\}$$ and $${\rm F}S(M_1,\ldots, M_d)=\left( S_s(M_1,\ldots, M_d)\right)_{s=1}^\infty.$$ Here, $\otimes$ denotes the Kronecker-/tensorproduct of matrices. The set $S_1(M_1,\ldots, M_d)$ is known as a classical {\it spectrahedron}. Such sets are precisely the feasible sets of semidefinite programming. The free spectrahedron  ${\rm F}S(M_1,\ldots, M_d)$ is a non-commutative extension, precisely in the spirit as above. For free spectrahedra, there exists a nice Positivstellensatz. As in Theorems \ref{hth} and \ref{box} above, we see that the natural notion of positivity in the non-commutative setup is strong enough to provide the best possible algebraic 
certificate (we do not cite the most general result and suppress some minor technical details for better readability):

\begin{Thm}[\cite{he2}] Let $M_1,\ldots, M_d\in {\rm Her}_r(\C)$ and $p\in \C\langle z_1,\ldots, z_d\rangle_h$. If $$p(\underline A)\geqslant 0$$ for all $\underline A\in S_s(M_1,\ldots, M_d)$ and all $s\geq 1$, in other words if $p$ is nonnegative on  the free spectrahedron ${\rm F}S(M_1,\ldots, M_d)$, then there exists a representation $$p=\sum_i q_i^*q_i+ \sum_j f_j^*Mf_j $$ where $  q_i\in \C\langle z_1,\ldots, z_d\rangle, f_j\in \C\langle z_1,\ldots, z_d\rangle^r$ and $$M:=z_1M_1+\cdots +z_dM_d\in {\rm Her}_r\left(\C\langle z_1,\ldots, z_d\rangle\right).$$
\end{Thm}

Sometimes facts about classical spectrahedra can only be learned by extending them to the non-commutative setup. One such instance is the  {\it containment problem} for spectrahedra, a problem appearing in different areas of (applied) mathematics \cite{ktt}. Given $M_1,\ldots, M_d\in {\rm Her}_r(\C)$ and $N_1,\ldots, N_d\in {\rm Her}_t(\C)$, how can one check efficiently  whether \begin{equation}\tag{A}S_1(M_1,\ldots, M_d)\subseteq S_1(N_1,\ldots, N_d)\end{equation} holds? Since spectrahedra are generalizations of polyhedra (which appear in the case of commuting coefficient matrices), this includes the problem of polyhedral containment. An important algorithm to solve this problem was proposed in \cite{bt}. Instead of checking (A) once checks \begin{equation}\tag{B} \exists V_1,\ldots, V_n\in {\rm Mat}_{r,t}(\C)\colon\quad \sum_jV_j^*M_iV_j=N_i \mbox{ for } i=1,\ldots, d.\end{equation} It is obvious that (B) implies (A).  Condition (B) can be transformed into a semidefinite optimization problem, and thus often solved efficiently. It was however known that (A) and (B) are not equivalent, so the answer to (B) could be {\it no}, where the answer to (A) is {\it yes}. A much better understanding of the method was gained through the following result (again we suppress some minor technical details):

\begin{Thm}[\cite{hkinf}]\label{hkkk} Condition (B) is equivalent to \begin{equation}\tag{A'}{\rm FS}(M_1,\ldots, M_d)\subseteq {\rm FS}(N_1,\ldots, N_d).\end{equation} Inclusion is meant level-wise here, i.e. $S_s(M_1,\ldots, M_d)\subseteq S_s(N_1,\ldots, N_d)$ for all $s\geq 1$.
\end{Thm}

This result mostly relies on Choi's characterization of completely positive maps between matrix algebras \cite{choi}. The insight of Theorem \ref{hkkk} can now be used to determine instances in which (A) and (B) are equivalent nonetheless. For this let $C\subseteq \R^d$ be a convex cone. There is one smallest and one largest matrix convex set with $C$ at level one. Indeed define $$C^{\min}_s:=\left\{ \sum_i c_i^t\otimes P_i\mid c_i\in C, P_i\in {\rm Her}_s(\C), P_i\geqslant 0\right\}$$ and $$C^{\max}_s:=\left\{ \underline A\in {\rm Her}_s(\C)^d\mid v^*\underline Av\in C\ \mbox{ for all } v\in \C^s\right\}.$$ Then $$C^{\min}:=\left(C^{\min}_s\right)_{s=1}^\infty\quad \mbox{and}\quad C^{\max}:=\left(C^{\max}_s\right)_{s=1}^\infty$$ are easily checked to be the smallest/largest such matrix convex set. Now assume $$C=S_1(M_1,\ldots, M_d)\subseteq \R^d$$ is a (classical) spectrahedral cone with ${\rm F}S(M_1,\ldots, M_d)= C^{\min}$. In that case, condition (A) implies $${\rm F}S(M_1,\ldots, M_d)=C^{\min} \subseteq {\rm F}S(N_1,\ldots, N_d)$$  and thus (B), by Theorem \ref{hkkk}. On the other hand, if $$C^{\min}\subsetneq {\rm F}S(M_1,\ldots, M_d)$$ it can be shown by the non-commutative separation theorem from \cite{ew}, that there exist matrices $N_1,\ldots, N_d$ with $$C=S_1(M_1,\ldots, M_d)=S_1(N_1,\ldots, N_d)$$ and  $$ {\rm F}S(M_1,\ldots, M_d) \nsubseteq{\rm F}S(N_1,\ldots, N_d).$$ In such an instance the answer to (B) is no, whereas the answer to (A) is yes. So the method from \cite{bt} works reliably if any only if  ${\rm F}S(M_1,\ldots, M_d)$ is the smallest matrix convex cone over the classical spectrahedron $S_1(M_1,\ldots, M_d)$. Unfortunately, this happens very rarely, already for polyhedral cones:

\begin{Thm}[\cite{fnt}]\label{op}
Assume $C=S_1(M_1,\ldots, M_d)\subseteq\R^d$ is polyhedral. Then $${\rm F}S(M_1,\ldots, M_d)= C^{\min}$$ if and only if $C$ is a simplex cone, i.e. has only $d$ extremal rays.
\end{Thm}

The last theorem also has some surprising application in theoretical quantum physics. The state of a bipartite quantum system is usually described by a positive semidefinite matrix $$\rho\in{\rm Mat}_r(\C)\otimes{\rm Mat}_s(\C)\cong {\rm Mat}_{rs}(\C).$$ So $\rho$ can be written as $$0\leqslant\rho=\sum_{i=1}^n M_i\otimes A_i$$ with $M_i\in {\rm Mat}_r(\C), A_i\in{\rm Mat}_s(\C).$ Although $\rho$ is supposed to be positive semidefinite and in particular Hermitian, this is not necessarily true for the $M_i, A_i$. If there exists a different such representation where all the $M_i,A_i$ are positive semidefinite as well, then $\rho$ is called {\it separable}, otherwise it is {\it entangled}. The smallest possible $n$ in the representation of $\rho$ above is called the {\it tensor rank} of $\rho$. A corollary of Theorem \ref{op} now reads as follows:

\begin{Thm}[\cite{ddn}]Every bipartite quantum state of tensor rank $2$ is separable.
\end{Thm}
 In fact $$0\leqslant\rho=M_1\otimes A_1+M_2\otimes A_2$$ just means that $(A_1,A_2)\in S_s(M_1,M_2)$. Now since the  convex cone $$C=S_1(M_1,M_2)\subseteq \R^2$$ is automatically a simplex cone, we obtain $(A_1,A_2)\in C_s^{\min}$ from Theorem \ref{op}. Writing down a representation in this smallest matrix convex cone and using bilinearity of the tensor product immediately implies the result.

Let us close with a result about non-commutative polytopes and polyhedra. The theorem of Minkowski-Weyl (see for example \cite{schr}) states that every polyhedral cone $C\subseteq \R^d$ is finitely generated, and vice versa. In other words, the notions {\it polyhedral} and {\it polytopal} coincide for convex cones. Now a short contemplation reveals that  $C^{\min}$ is a good generalization of the notion  {\it polytope/finitely generated} to the non-commutative setup, whereas $C^{\max}$  corresponds to the {\it polyhedral} notion. Interestingly, these two notions differ almost always, already at the first level of non-commutativity:

\begin{Thm}[\cite{fnt,sh,hn}] Let $C\subseteq \R^d$ be a convex cone.

($i$) If $C$ is a simplex cone, then $C^{\min}=C^{\max}$. Otherwise $C^{\min}\neq C^{\max}$.

($ii$)  If $C$ is polyhedral but not a simplex cone, then $C_2^{\min}\subsetneq C_2^{\max}.$
\end{Thm}

As a concluding remark, we note that the methods used in the non-commutative setup differ quite strongly from the ones in the commutative theory. Many of the results are proven by functional-analytic methods, such as GNS-constructions, dilations, and the theory of completely positive maps and operator algebras. Sometimes results and examples from group theory and the theory of $C^*$-algebras can be useful. All in all, the whole area is not yet mature, many interesting results and methods are hopefully developed in the coming years.

\end{document}